\documentclass[%
 aip,
 amsmath,amssymb,
reprint,%
,cha
]{revtex4-1}

\usepackage{graphicx}
\usepackage{dcolumn}
\usepackage{bm}

\usepackage[utf8]{inputenc}
\usepackage[T1]{fontenc}
\usepackage{mathptmx}
\usepackage{xcolor}
\graphicspath{{./figures/}}
\usepackage{tabularx} 
\usepackage{booktabs} 
\usepackage{subcaption} 
\usepackage{bbm}

\usepackage{listings}

\definecolor{codegreen}{rgb}{0,0.6,0}
\definecolor{codegray}{rgb}{0.5,0.5,0.5}
\definecolor{codepurple}{rgb}{0.58,0,0.82}
\definecolor{backcolour}{rgb}{0.95,0.95,0.98}

\lstdefinestyle{mystyle}{
    backgroundcolor=\color{backcolour},   
    commentstyle=\color{codegreen},
    keywordstyle=\color{magenta},
    numberstyle=\tiny\color{codegray},
    stringstyle=\color{codepurple},
    basicstyle=\ttfamily\footnotesize,
    breakatwhitespace=false,         
    breaklines=true,                 
    captionpos=b,                    
    keepspaces=true,                 
    showspaces=false,                
    showstringspaces=false,
    showtabs=false,                  
    tabsize=2
}

\def\cvdots {\omit\span\omit \hfil$\vdots$\hfil} 

\begin{document}

\preprint{AIP/123-QED}

\title{Limitations of the Generalized Pareto Distribution-based estimators for the local dimension}

\author{Ignacio del Amo}
\email{i.b.del-amo@exeter.ac.uk}
\affiliation{Department of Mathematics and Statistics, University of Exeter, EX4 4QF, United Kingdom}

\author{George Datseris}
\affiliation{Department of Mathematics and Statistics, University of Exeter, EX4 4QF, United Kingdom}

\author{Mark Holland}
\affiliation{Department of Mathematics and Statistics, University of Exeter, EX4 4QF, United Kingdom}

\date{\today}


\begin{abstract}
Two dynamical indicators, the local dimension and the extremal index, used to quantify persistence in phase space have been developed and applied to different data across various disciplines. These are computed using the asymptotic limit of exceedances over a threshold, which turns to be a Generalized Pareto Distribution in many cases. However the derivation of the asymptotic distribution requires mathematical properties which are not present even in highly idealized dynamical systems, and unlikely to be present in real data. Here we examine in detail issues that arise when estimating these quantities for some known dynamical systems with a particular focus on how the geometry of an invariant set can affect the regularly varying properties of the invariant measure. We demonstrate that singular measures supported on sets of non-integer dimension are typically not regularly varying and that the absence of regular variation makes the estimates resolution dependent. We show as well that the most common extremal index estimation method is ambiguous for continuous time processes sampled at fixed time steps, which is an underlying assumption in its application to data.

\end{abstract}

\maketitle


\section{Introduction}

In recent years there has been a growing interest in Extreme Value Theory based methods to obtain estimates of dynamical quantities. These methods stem from developments in the analysis of extreme events and its relationships with dynamical systems. Same as the sum of random variables has an asymptotic distributional limit described by the Central Limit Theorem, less known theorems describe the distributional limits of extreme events, allowing the quantification of probabilities that relate the magnitude of events with their frequency. These results have been extended to cover extremes that arise in observables computed along trajectories of chaotic dynamical systems. By defining an observable that depends on the distance between a reference point and the trajectory of the system in state space, one can obtain a representation of the extreme events as close recurrences to the reference point, which relate the dynamical and ergodic properties of the system to the statistics of the extremes. These methods have been extensively used for atmospheric data \cite{alberti2023dynamical,faranda2017dynamical,faranda2017dynamicalB,faranda2020diagnosing,hochman2019new,messori2017dynamical,rodrigues2018dynamical}, and fractal dimension estimation\cite{datseris2023estimating}.

There exist mainly two kinds of Extreme Value Laws. The asymptotic limit of the maximum value in an observed data series partitioned in blocks follows a Generalized Extreme Value distribution (GEV), meanwhile the asymptotic limit of the exceedances over a high threshold follow a Generalized Pareto Distribution (GPD). These laws have been shown to apply also to apply to extremes of observables along orbits of dynamical systems under sufficiently fast decay of correlations\cite{freitas2012extremal,freitas2015speed}.
The theory for estimating the local dimension through extreme value theory has been developed theoretically mostly for the GEV case \cite{faranda2012generalized,lucarini2012extreme,faranda2011numerical,faranda2014extreme}, with only a couple of papers dealing with the GPD estimator \cite{lucarini2014towards,lucarini2012universal}. However, as typically the GPD approach requires less data, it is the one that has been mostly used in applications. But this approach requires several abstract mathematical properties to hold, which are different from the ones required by the GEV approach. On one hand, there are some statistical properties of the dynamical system such as ergodicity and stationarity. On the other hand there are some geometrical conditions on the invariant measure of the system, such as the existence of the local dimension and regular variation.


In this paper we focus in numerically investigating the regular variation property, which is a necessary condition, but it is only mentioned sometimes, and to our knowledge, has never been checked for any system or dataset in the literature dealing with the GPD method. Such is the case that even some texts that explicitly intend to illustrate the use of this method, do so in systems which are not regularly varying. 
For measures defined in invariant sets of dynamical systems it is complicated to obtain an analytical result regarding its regular variation properties.
Sometimes an estimation can be done for the measures supported on self similar attractors analytically, but still this only leads to a case by case check. Here we show how the regular variation of the measure relates to the geometry of the attractor, and how the lack of it makes the estimates for the local dimension scale dependent for many commonly used dynamical systems, putting special emphasis in those which have been previously used in the literature.


\section{Background}

The theory behind the algorithm to compute local dimensions based in the GPD \cite{lucarini2012universal,lucarini2014towards} is the result of applying the Peaks over Threshold (PoT) method to a stochastic process of the form
\begin{align}\label{eq:hittingprocess}
    X_i = -\log(||T^i\mathbf{x}_0-\zeta||)
\end{align}
where $T$ is a map generally representing the forwards time evolution of a dynamical system for some time $\delta t$, $\textbf{x}_0$ and $\zeta$ are points in phase space belonging to an invariant set where $\zeta$ is the reference point and has to be chosen independently from $\mathbf{x}_0$ and $\{T^i\mathbf{x}_0\}_{i=0}^N$ is an orbit. From now on we will refer to this algorithm as Exceedance-Based Dimension (EBD) algorithm.
The PoT method consists on taking a high threshold and defining as extreme events the exceedances over the threshold, i.e. the values of Eq.~\eqref{eq:hittingprocess} that surpass the threshold $g_q$. This notation for the threshold is meant to emphasize that it can be defined via a quantile, and we use the convention that if $q = 0.99$ then the threshold $g_q$ leaves $99\%$ of data below and only $1\%$ above. Then to find the probabilities corresponding to excesses over the threshold means that we are interested in computing the asymptotic excess distribution
\begin{align}
    H(u) = \lim_{g_q\to\infty}\mathbb{P}(X_n-g_q\leq u|u>g_q).
\end{align}
Note that since the process $X_i$ is defined as a decreasing monotonic function of the distance between the orbit and a point, meaning that a high value of the process corresponds to a close "recurrence" to a ball near the reference point. Under the assumptions of ergodicity and stationarity, there is a invariant measure $\mu$ that describes the probability of finding the process in a particular region of state space, and thus the probability of a exceedance over the threshold is $\mu(B_{e^{-g_q}}(\zeta))$, where we use the notation $B_r(\mathbf{x})$ to denote a ball of radius $r$ centered in the point in phase space $\mathbf{x}$. The excess distribution function becomes
\begin{align}
    \mathbb{P}(X_n-g_q\leq u|u>g_q)= 1 - \frac{\mu(B_{e^{-u-g_q}}(\zeta))}{\mu(B_{e^{-g_q}}(\zeta))}
\end{align}
since is the difference in volume between two balls of radii $e^{-g_p}$ and $e^{-u-g_p}$, normalized by the volume of a ball of radius $e^{-g_p}$.

The local dimension is a dynamical quantity that quantifies how the measure of a ball scales with its radius. It is defined as 
\begin{align}\label{eq:localdimensiondefinition}
    \Delta^l_\zeta = \lim_{r\to 0}\frac{\log \mu(B_r(\zeta))}{\log r}
\end{align}
for the points in the support of $\mu$ for which the limit exists. 
If we assume that the measure is non-atomical and the local dimension exists almost everywhere, then the measure of a ball of radius $r$ is a smooth function of the radius. In particular, for small radius  $\mu(B_r(\zeta))\simeq l(r) r^{\Delta^l_\zeta}$ where $l(r)$ is a function of the radius such that $\lim_{r\to0} \log f(r)/\log r = 0$. To conclude the derivation now we need to assume that $\mu$ is a regularly varying measure, meaning that it has the property
\begin{align}\label{eq:regularindex}
    \lim_{r\to 0}\frac{\mu(B_{br}(\zeta))}{\mu(B_r(\zeta))} \to b^\gamma
\end{align}
for any $0<b\leq 1$, where $\gamma$ is called the regularly varying index. Note that if the local dimension exists, we have that
\begin{align}
    \frac{\mu(B_{br}(\zeta))}{\mu(B_r(\zeta))} \simeq \frac{l(br)}{l(r) } \frac{(br)^{\Delta^l_\zeta}}{r^{\Delta^l_\zeta}}  \to b^{\Delta^l_\zeta}
\end{align}
and implies that the function $l(r)$ is slowly varying, i.e.  $\lim_{r\to 0} l(br)/l(r)\to 1$ for any $0<b\leq 1$. Thus, the index $\gamma$ is equal to the local dimension around the point $\zeta$. Thus the local dimension and the index of regular variation are intrinsically connected.
With the slow variation property, the limit in the excess distribution function can be solved, and reduces to the distribution function of an exponential distribution with parameter $1/\Delta^l_\zeta$
\begin{align}
\begin{aligned}
    H(u) &= 1 - \lim_{g_p\to \infty} \frac{\mu(B_{e^{-u-g_p}}(\zeta))}{\mu(B_{e^{-g_p}}(\zeta))}\\
 & = 1 - \lim_{g_p\to \infty} \frac{l(e^{-g_p}e^{-u})(e^{-u-g_p})^{\Delta^l_\zeta}}{l(e^{-g_p})(e^{-g_p})^{\Delta^l_\zeta}}\\
    &= 1-e^{-\Delta^l_\zeta u}.
\end{aligned}
\end{align}
This gives us a powerful way of estimating the local dimension. On top of that, for systems which are exact dimensional, meaning that the local dimension is constant $\mu$-almost everywhere, the information dimension is equal to the value that the local dimension takes $\mu$-almost everywhere\cite{ott2002chaos}. In summary, the EBD algorithm requires ergodicity, stationarity, existence of the local dimension $\mu$-almost everywhere and regular variation of the measure to work.

In experimental applications one may only have one realization of the data series, which makes choosing an independent reference point difficult or wasteful, since one may have to discard many data points to decorrelate a reference point to the rest of the orbit. In practice, what is done in applications \cite{faranda2017dynamical} is to use as reference point a point in the orbit $\zeta = T^j x_0$ for some $j$, and thus one actually uses processes of the form 
\begin{align}\label{eq:algorithm1process}
    Y_i^j = -\log(||T^i\boldsymbol{x}_0-T^j\boldsymbol{x}_0||).
\end{align}
under the assumption that their extremes asymptotically behave as the extremes of the process given by Eq. \eqref{eq:hittingprocess}. Note that the reference point is not uncorrelated anymore to the points in the trajectory, and this makes it a non stationary process, since correlations to the reference point exist but decay over time so the process is asymptotically stationary. For fast mixing processes the orbit decorrelates quickly from the initial condition, and hence if $j\ll i$ we recover the same asymptotic result with $T^j\textbf{x}_0$ replaced by $\zeta$. It is unclear however what happens when $i\ll j$, and the reference point is in the distant future of the orbit instead of the distant past. The discussion of the effects of substituting process \eqref{eq:hittingprocess} by process \eqref{eq:algorithm1process} is left for future work.
In applications \cite{faranda2017dynamical,alberti2023dynamical,faranda2017dynamicalB,faranda2020diagnosing,hochman2019new,rodrigues2018dynamical,messori2017dynamical}, the output of the EBD algorithm is usually paired with an estimation of a quantity called the extremal index, usually denoted by $\theta$. The extremal index arises in the Generalized Extreme Value distributions in the presence of short range dependencies and can be seen as a number that describes the inverse of the mean size of clusters of extremes \cite{moloney2019overview}. Note that the derivation of the exponential law for the exceedances does not need a particular dependency structure, and does not feature the extremal index in the limiting distribution. However in some papers \cite{faranda2017dynamical,alberti2023dynamical,faranda2017dynamicalB,messori2017dynamical,faranda2020diagnosing,hochman2019new} discussing applications it is included in the distribution function without a mathematical justification. The PoT approach erases the dependencies by considering the points over threshold without regard for its position in the process and thus shuffling the time series randomly will produce the exact same extremes. The extremal index appears naturally in the Block Maxima approach since the clusters typically fall within a block, and thus only one extreme is chosen from each cluster. 

The following example shows that the Peaks over Threshold approach does not require the extremal index. Consider a process $V_i$ defined by i.i.d. draws from an exponential distribution of parameter $\lambda$ and $U_i = \max\{V_{i-1},V_i\}$. The distribution function of $V_i$ is $F_V(y)= 1-e^{-y/\lambda}$, and for $U_i$ is $F_U(y) = F_V(y)^2$.
Then the asymptotic distribution of the maxima of $V_i$ is 
\begin{align*}
    \lim_{n\to\infty}\mathbb{P}\left(\max\{V_1,\dots,V_n\}<(y+\log n)/\lambda\right) 
    & = e^{-e^{-y}}
\end{align*}
whereas for $U_i$ is 
\begin{align*}
    \lim_{n\to\infty}\mathbb{P}\left(\max\{U_1,\dots,U_n\}<(y+\log n)/\lambda\right) 
    & = e^{-\theta e^{-y}}
\end{align*}
where $\theta=1/2$ is the extremal index. The extremal index appears due to the short time correlations introduced by the definition of $U_i$, and it is equal to $1/2$ particularly because of $F_U(y) = F_V(y)^2$. For the Peaks over Threshold approach however

\begin{align*}
\begin{aligned}
    \lim_{g_q\to\infty}\mathbb{P}(V_n-g_q\leq y|y>g_q)& = \lim_{g_q\to\infty}\frac{F_V(y-g_p)-F_V(g_p)}{1-F_V(g_p)}\\
    &= 1-e^{-x/\lambda}
\end{aligned}
\end{align*}
\begin{align*}
\begin{aligned}
    \lim_{g_q\to\infty}\mathbb{P}(U_n-g_q\leq y|y>g_q)& = \lim_{g_q\to\infty}\frac{F_V(y-g_p)^2-F_V(g_p)^2}{1-F_V(g_p)^2}\\
    &= 1-e^{-x/\lambda}.
\end{aligned}
\end{align*}

While the extremal index can be seen as statistical quantity that is intrinsically linked to the process $U_i$ and can be estimated from data, it does not appear in the asymptotic distribution of extremes through the PoT method.

One further comment that should be made is that typically the set of reference points in a chaotic invariant set that give raise to an extremal index different than 1 is the set of stable and periodic points\cite{lucarini2016extremes}, and thus a set of $\mu$ measure zero, unlikely to be observed in most applications.

\section{Discrete Systems}
\subsection{Cantor Shift}

We start our exposition with the Cantor Shift map for two reasons. One is that is such a simple system that it can be shown that its invariant measure is not regularly varying. The second one is that it has been used in the literature \cite{lucarini2016extremes} to illustrate the EBD method, recognizing difficulties to apply it but without identifying it as a system with a non regularly varying invariant measure.

The middle third Cantor set is usually built through a Iterated Function System construction. A function system $f_C$ that maps the unit interval $C_0 =[0,1]$ to the union $C_1=[0,1/3]\cup[2/3,1]$ is iterated repeatedly, and then the middle third cantor set is defined as
\begin{align*}
    & C_n = f_C^n([0,1])\\
    & C_\infty = \bigcap_{i=0}^\infty C_i,
\end{align*}
for a more detailed construction and explanation, see \cite{falconer2007fractal}.

A way to build a dynamical system in the middle third Cantor set $C_\infty$ is given by symbolic dynamics. We identify a point in the $C_\infty$ with a semi-infinite sequence in the space $\{0,2\}^\mathbb{N}$. A point in this space can be regarded as an decimal expansion $\omega= 0.a_1a_2a_3...$ where $a_i\in\{0,2\}$. Then the shift map acting on $\omega$ is defined as the map that deletes the first element of the expansion and moves the decimal point to the right, this is $T_c(0.a_1a_2a_3...)= 0.a_2a_3a_4...$. The shift can be related to the middle third Cantor Set via the embedding
\begin{align}
    P(0.a_1a_2a_3...) = \sum_{i=1}^\infty \frac{a_i}{3},
\end{align}
and this implies that the shift can be seen as a dynamical system supported in the Cantor set, which is an invariant set for this transformation. If $\omega$ is chosen such that the asymptotic density of $0$'s and $2$'s is $1/2$, for example generating it with a Bernoulli process, then the $(1/2,1/2)$ Bernoulli measure $\mu_C$ is the relevant invariant measure for the dynamics. This measure is akin to a uniform measure but that has $C_\infty$ as its support, which makes it singular with respect Lebesgue measure, given that the Lebesgue measure of the Cantor set is zero and the Bernoulli measure of $[0,1]\setminus C_\infty$ is zero as well.

\begin{figure}
    \centering
    \includegraphics[width = \columnwidth]{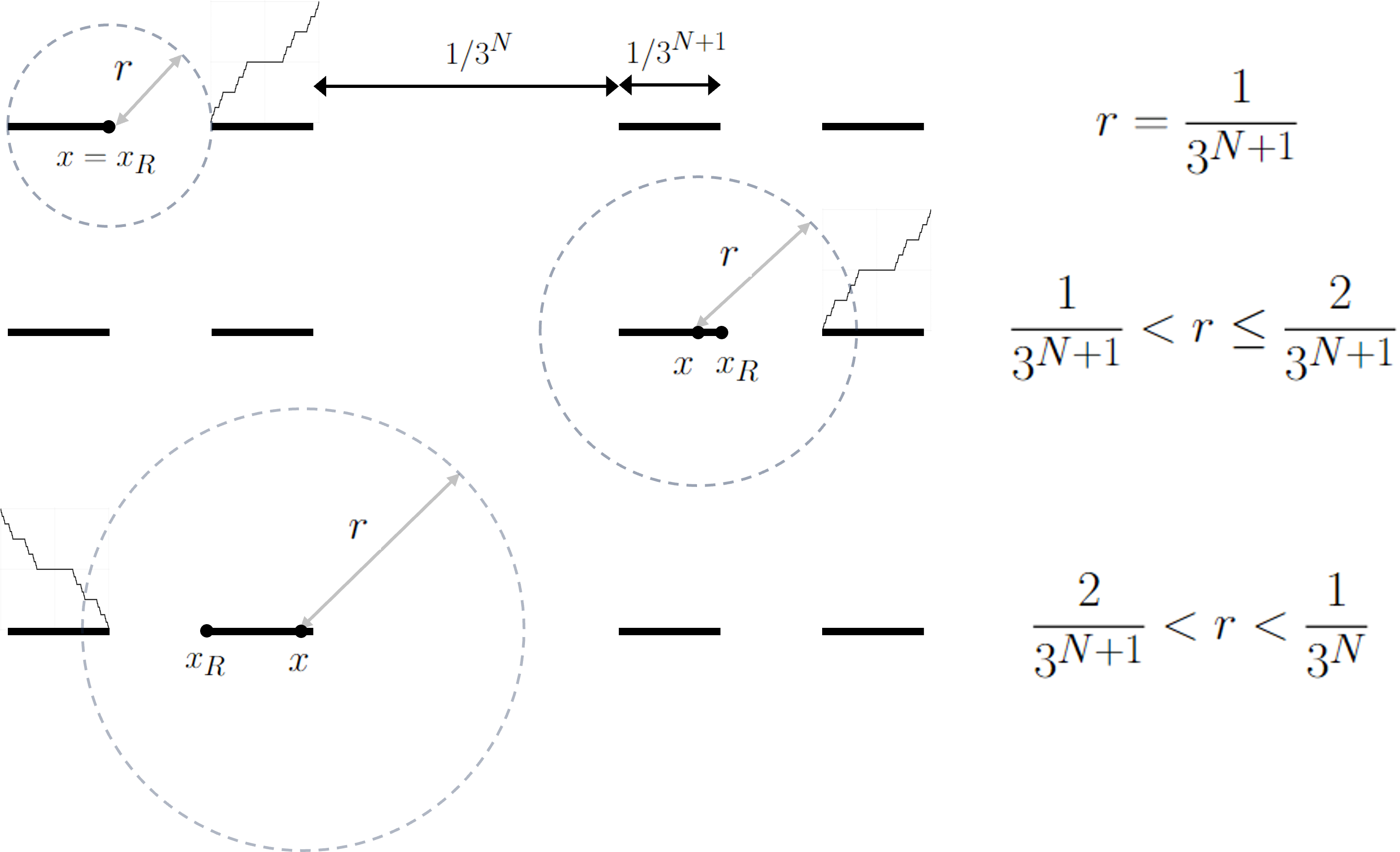}
    \caption{Visual aid to understand the computation of the measure of a ball with the (1/2,1/2) Bernoulli measure. The segments represent a subset of set $C_{N+1}$, since one of the segments is completely inside the ball always, further detail is not needed, however, when the radius is not a power of 1/3, the intersection with a nearby segment makes infinite detail necessary. This takes the form of the Cantor ternary function $C(\cdot)$, represented on top of the pertinent segments.(Change $x$ for $\zeta$ in figure)}
    \label{fig:BernoulliMeasure}
\end{figure}

The measure of ball a radius $r < 1$ and centered in a point $\zeta\in C_\infty$ can be computed in the following way. Note that for some $N\in \mathbb{N}$ we have that $1/3^{N+1}\leq r<1/3^N$. Regardless of the position of the center we have that $1/2^N\geq\mu_b(B_r(\zeta)) \geq 1/2^{N+1}$ because the ball will always contain a segment of $C_{N+1}$, with equality holding for all $x$ if $r$ is exactly  $1/3^{N+1}$. 
When $r\leq 2/3^{N+1}$, the measure of the ball will depend on the position of the point, since there will be a subset of positive measure of $C_\infty$ which is at $1/3^{N+1}$ distance, but only in one of the sides, and there might be an intersection. Let $x_e$ be the endpoint of the subset $C_{N+1}$ such that $\zeta\in C_{N+1}$ which is at distance $1/3^{N+1}$ of the closest subset of $C_{N+1}$. Then the measure of the ball is $1/2^{N+1}(1+C(3^{N+1}(r-1/3^{N+1}-|\zeta-x_e|))$ where $C(x)$ is the Cantor ternary function (also called Devil's staircase) whenever $r-1/3^{N+1}-|\zeta-x_e|$ is positive, or 0 otherwise.
When $r\geq 2/3^{N+1}$ the ball will always contain some positive measure subset of the nearest subset. The measure of the ball as before can be expressed again as $1/2^{N+1}(1+C(3^{N+1}(r-1/3^{N+1}-|\zeta-x_e|))$ with $C(x)+1$ when the $x\geq1$.

The Bernoulli measure is not regularly varying. Note that because of the gaps in the support of the measure, there are values of $b$ for which the measure does not change value at all. Let $b = 1/2$ and $\zeta - x_e = 1/3^{N+1}$ meaning $\zeta$ is the opposite endpoint. Note that the endpoints will always belong to $C_\infty$ by construction. If $r=2/3^{N+1}$ then Eq.\eqref{eq:regularindex} can only be true if the index $\gamma=0$. Values of $b$ and $r$ such that this holds can be found for any point $\zeta\in C_\infty$, so decreasing sequences of $r$ such that this holds can be found. 

By taking any point $\zeta$, any radius and $b=1/3$, then due to self-similarity we have

\begin{align*}
    \frac{\mu(B_{br}(\zeta))}{\mu(B_r(\zeta)} &= \frac{1/2^{N+2}(1+C(3^{N+2}(br-1/3^{N+2}-|\zeta-x_{R_b}|))}{1/2^{N+1}(1+C(3^{N+1}(r-1/3^{N+1}-|\zeta-x_R|))}\\
    &= 1/2 = b^{\frac{\log 2}{\log 3}},
\end{align*}
thus the limit in Eq.~\eqref{eq:regularindex} does not exist. Note that using instead other $(p,1-p)$-Bernoulli measure, or other self similar 1D zero measure Cantor set does not change this fact.

\begin{figure}
    \centering
    \includegraphics[width = \columnwidth]{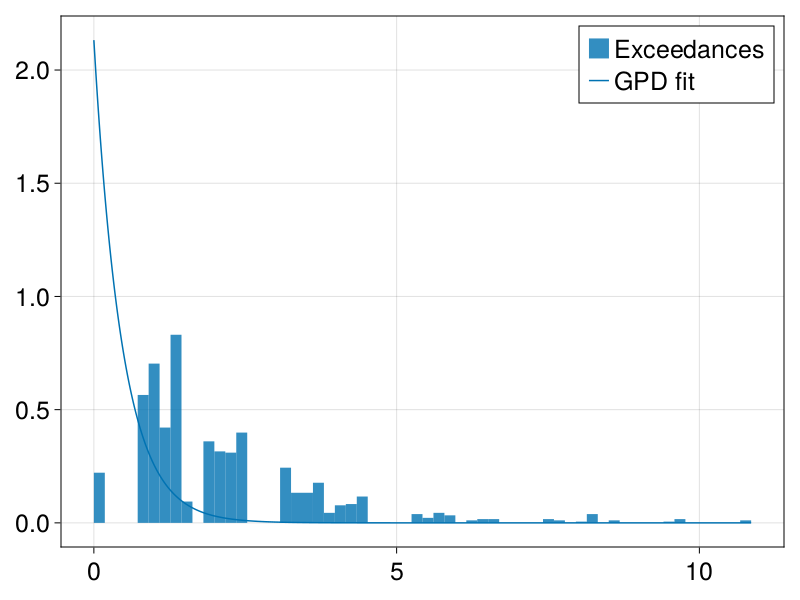}
    \caption{Distribution of exceedances for the Cantor shift map, and a GPD distribution fitted to it by assuming that is an exponential distribution with parameter equal to the inverse of the average of the exceedances.}
    \label{fig:CantorGPD}
\end{figure}

Figure~\ref{fig:CantorGPD} shows the result of taking the exceedances from the Cantor shift map and fitting a GPD distribution to them according to the algorithm. The holes in the measure produce in turn holes in the distribution of exceedances, as certain values are not possible.
 \begin{figure*}
         \centering
     \includegraphics[height=0.9\textheight]{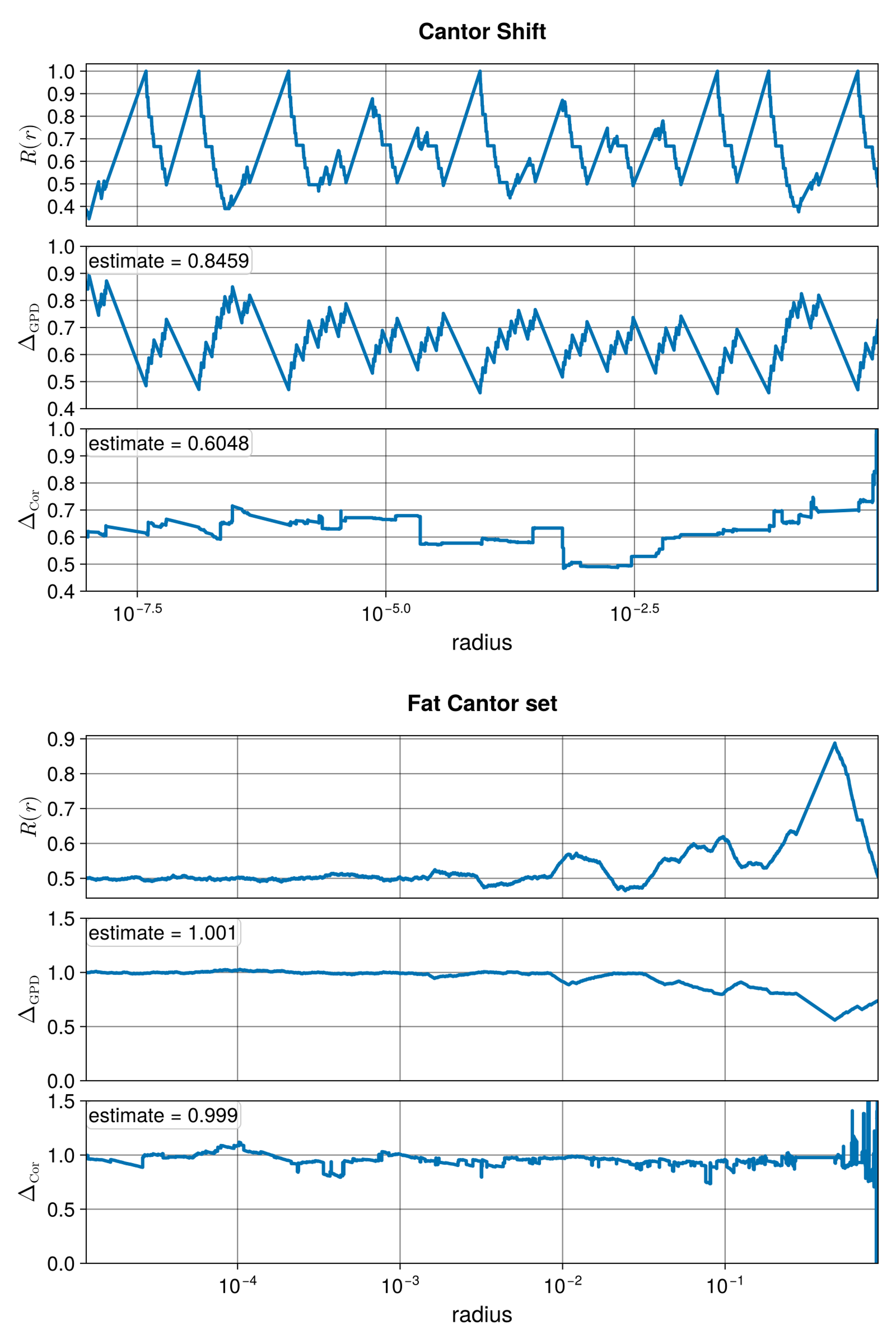}
     \caption{The top panel is formed by the first three graphs and shows the quantity $R(r)$ in the top, the estimation of the local dimension via the EBD method and the estimation of the local dimension via the correlation method. The bottom panel is formed for the corresponding graphs for the fat Cantor set of points which is invariant under the non-autonomous tent map.}
     \label{fig:Cantor_maps_variation}
 \end{figure*}

Figure~\ref{fig:Cantor_maps_variation} is prototypical of a kind of figure we will use repeatedly throughout the paper to analyse each system, and as such it deserves an in depth explanation. The top panel shows three graphs with linked axis, in the top graph a numerical computation of the quantity
\begin{align}\label{eq:regularly_varying_quantity}
    R(r) = \frac{\mu(B_{r/2}(\zeta))}{\mu(B_r(\zeta))} \simeq \frac{l(r/2)}{l(r)}\frac{(r/2)^{\Delta^l_\zeta}}{r^{\Delta^l_\zeta}}
\end{align}
for a broad range of decreasing radii is shown, which means that if the measure was regularly varying, we would see this quantity converge to the constant $2^{-\Delta_\zeta^l}$. The way in which this is computed is the following. A random reference point $\zeta$ and initial condition $\boldsymbol{x}_0$ in the invariant set are chosen, then the orbit is computed forward, keeping track of which are the 5000 highest exceedances, which correspond to the closest recurrences of the orbit to a ball centred on $\zeta$. When a higher exceedance is encountered, the smallest one is deleted, and so the radius of the ball decreases because the lowest exceedance corresponds to the point furthest away from $\zeta$. The middle panel shows the estimated local dimension around the reference point estimated with the points which lay in the ball through the EBD method described above. The bottom panel shows an estimation of the local dimension made through an approach based on correlations between points due to Grassberger-Procaccia\cite{grassberger1983characterization}. This estimator is given for comparison, since it is a method well established in the literature and does not depend on the regularly varying properties of the measure, see appendix \ref{ap:correlation_dimension_explanation} for more details.

The top panel of Figure~\ref{fig:Cantor_maps_variation} shows that the quantity described in \eqref{eq:regularly_varying_quantity} does not converge to $2^{-\log 2/\log3}\approx 0.646$. Instead it oscillates strongly for any of the orders of magnitude reached. Since the measure is not regularly varying, we see a non-vanishing factor that prevents the convergence of the asymptotic distribution to an exponential. The middle panel shows that the estimate of the local dimension is highly correlated with the concentration of the measure and oscillates as well, making the estimates of the local dimension resolution dependent. The correlation estimate for the dimension, although it also goes up and down a little, shows reasonable values and the final value is relatively close to the actual analytical solution $\log 2/\log3\approx 0.631$. 

\subsection{Fat Cantor set}

An interesting edge case is given by the so called ``fat fractals''. These are sets of positive Lebesgue measure and thus of integer dimension, but that might have a very intricate fine structure similar to a fractal set. The example that we show here is a Cantor set, a nowhere dense set of isolated points and thus topologically equivalent to the middle third Cantor set, but it has positive Lebesgue measure. The set is produced via a dynamic construction, using a non-autonomous modification of the tent map.

\begin{align*}
    f_C^{(n)}(x) = 
    \begin{cases}
        2\left(1+2^{-n-1}\right)x & \text{ if $x\leq1/2$}\\
        2\left(1+2^{-n-1}\right)(-x+1) & \text{ if $1/2<x$}
    \end{cases}
\end{align*}

The value of $f_C^{(n)}$ exceeds 1 for the points situated in the interval $I_n =((2(1+2^{-n-1}))^{-1},1-(2(1+2^{-n-1}))^{-1})$. The points that are mapped outside of the interval $[0,1]$ escape to infinity, thus in the first iteration the points in $I_1$ escape, in the second iteration the points in ${\big(}f_C^{(1)}{\big)}^{-1}(I_2)$ escape, and then the invariant set of points which do not leave the interval as $n\to\infty$ can be written as

\begin{align*}
    A_C = ([0,1]\setminus I_1) \setminus \bigcup_{n=2}^\infty {\big(}f_C^{(1)}{\big)}^{-1}\circ ... \circ {\big(}f_C^{(n-1)}{\big)}^{-1} (I_{n}).
\end{align*}

In the $n$th iteration the Lebesgue measure of the remaining points diminishes by $2^{n-1}$ sets of measure $Leb(I_n)=1/(2^{n+1}+1)$ contracted by the inverse dynamics by a factor of $\prod_{i=1}^{n-1} 1/(2(1+2^{-i-1}))$. The Lebesgue measure of the invariant set can be estimated 
\begin{align*}
    Leb(A_C) & = 1 - 1/5 - \sum_{n=2}^\infty \frac{2^{n-1}}{2^{n+1}+1}\prod_{i=1}^{n-1} \frac{1}{2(1+2^{-i-1})}\\
    & > 4/5 - \sum_{n=2}^\infty \frac{2^{n-1}}{2^{n+1}}\prod_{i=1}^{n-1} \frac{1}{2}\\
    & = 4/5-1/4
\end{align*}
and thus is bounded away from zero. 

We generate the points in this set by drawing 5 initial conditions with a uniform distribution in the unit interval and iterating them forward 20 times, then taking points that remain in the interval and iterating them back 20 times again.
Since the set in this case is not generated by an orbit, there is no natural measure sitting atop this set, but since the EBD method utilizes only the geometrical information and not the order of the data we can randomly shuffle the points in the set and still recover the same answer. Note that the method we use to obtain the estimation of the regular variation properties might be sensitive to the order in which the data is distributed, but for a discrete chaotic system we expect that the time between recurrences to a small ball greatly exceeds the decorrelation time of the system, except in some exceptional cases like periodic or indifferent points, and thus the recurrences can be taken as independent from each other. The measure we effectively use is again the $(1/2,1/2)$ Bernoulli measure. Starting with a uniform distribution and because the maps are symmetrical respect to $x=1/2$, each iteration removes the same amount of measure from each side, and the measure of each $n$-cylinder set (corresponding to a $(i/2^n,(i+1)/2^n)$ interval) is $2^{-n}$, same as for the middle third Cantor set described above.


Despite the holes in the set and the topological equivalence to the Cantor Set, the $(1/2,1/2)$ Bernoulli measure here is regularly varying. This can be seen by the convergence of the quantity $R(r)$ (\eqref{eq:regularly_varying_quantity}) in the top graph of the bottom panel of figure~\ref{fig:Cantor_maps_variation} to the expected value of $2^{-1}$. The estimates of the local dimension also seem to converge to the right answer, which is 1. 
Due to the different embedding from the space $\{0,2\}^\mathbb{N}$ into $[0,1]$, the measure theoretic properties of the set change drastically, allowing it to have positive Lebesgue measure and making it regularly varying.

\subsection{H\'{e}non Map}

The Cantor shift map serves as an example of the consequences of the lack of regular variation. Although it is a very abstract and artificial example, helps to shed light on a widespread problem in hyperbolic systems, which often have singular measures whose support lies on a zero Lebesgue measure set. Cantor sets are also a models for the sections of attractors such as the H\'{e}non map.

\begin{figure*}
    \centering
    \includegraphics[width=\textwidth]{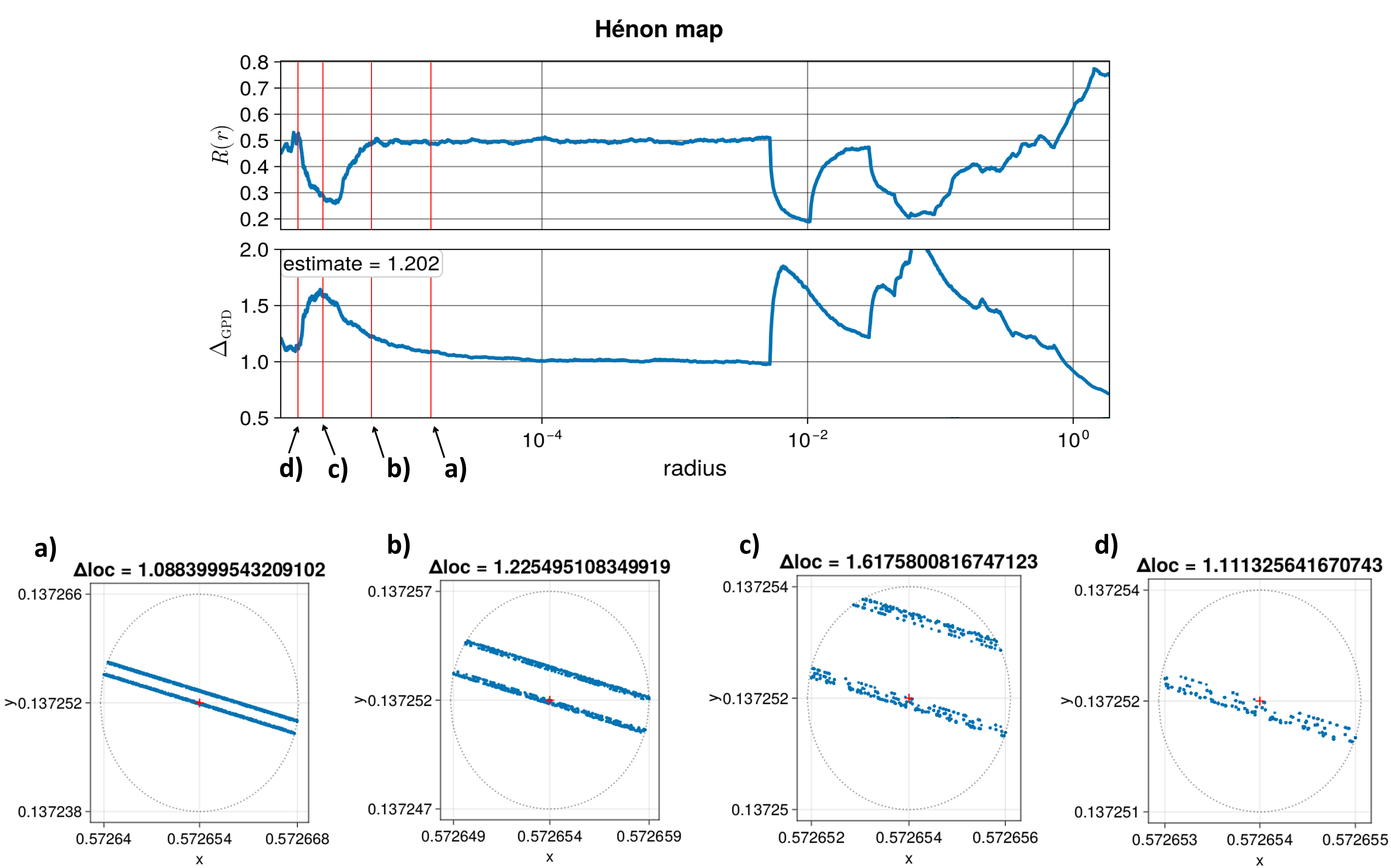}
    \caption{This figure shows the geometrical structure of the points that lay in the outer ball as the radius becomes small. The top graph shows the quantity $R(r)$, describing the ratio of the measure of two concentric balls, and bottom graph the corresponding estimation of the local dimension. Both plots seem to converge for a several orders of magnitude, but then an oscillation shows up. The vertical red lines indicate where in the zooming process we take the snapshots of the attractor shown below and labelled $a)$ to $d)$.
    }
    \label{fig:HenonZooms}
\end{figure*}

The H\'{e}non map is given by
\begin{align*}
    f_H(\textbf{x}) = (1-ax_1^2+x_2,bx_1)
\end{align*}
where we use the parameter values $a=1.4$, $b=0.3$, also called Benedicks-Carleson parameters. For this parameter values, the system has an attracting set, whose invariant measure is highly irregular, and it is hard to deal with analytically. Figure \ref{fig:HenonZooms} shows the same computation as in the previous figure but performed on a randomly selected point in the H\'{e}non attractor. On the right panels it also shows the points that lie inside the ball, the exceedances, and their geometrical disposition on four different stages of the zooming in process around $\zeta$. As the threshold becomes larger, the algorithm zooms more into the local structure of the attractor, and a branch that looks like a line from afar, doubles, becomes separate lines and eventually the one not containing $\zeta$ disappears and the one containing it starts dividing into more. This process causes the quotient of the measure of the balls and the estimation of the local dimension to oscillate strongly. Due to the structure of the attractor, we expect this problem to arise at all sufficiently small scales, for all points.

This implies in practice that the estimations obtained for the local dimension are scale dependent. Since the resolution of the data is given by the length of the trajectory, this is in most practical cases arbitrary, and thus we should be very careful when trying to extract information about an individual event based on the local dimension. 


\begin{figure*}
    \centering
    \includegraphics[height=0.9\textheight]{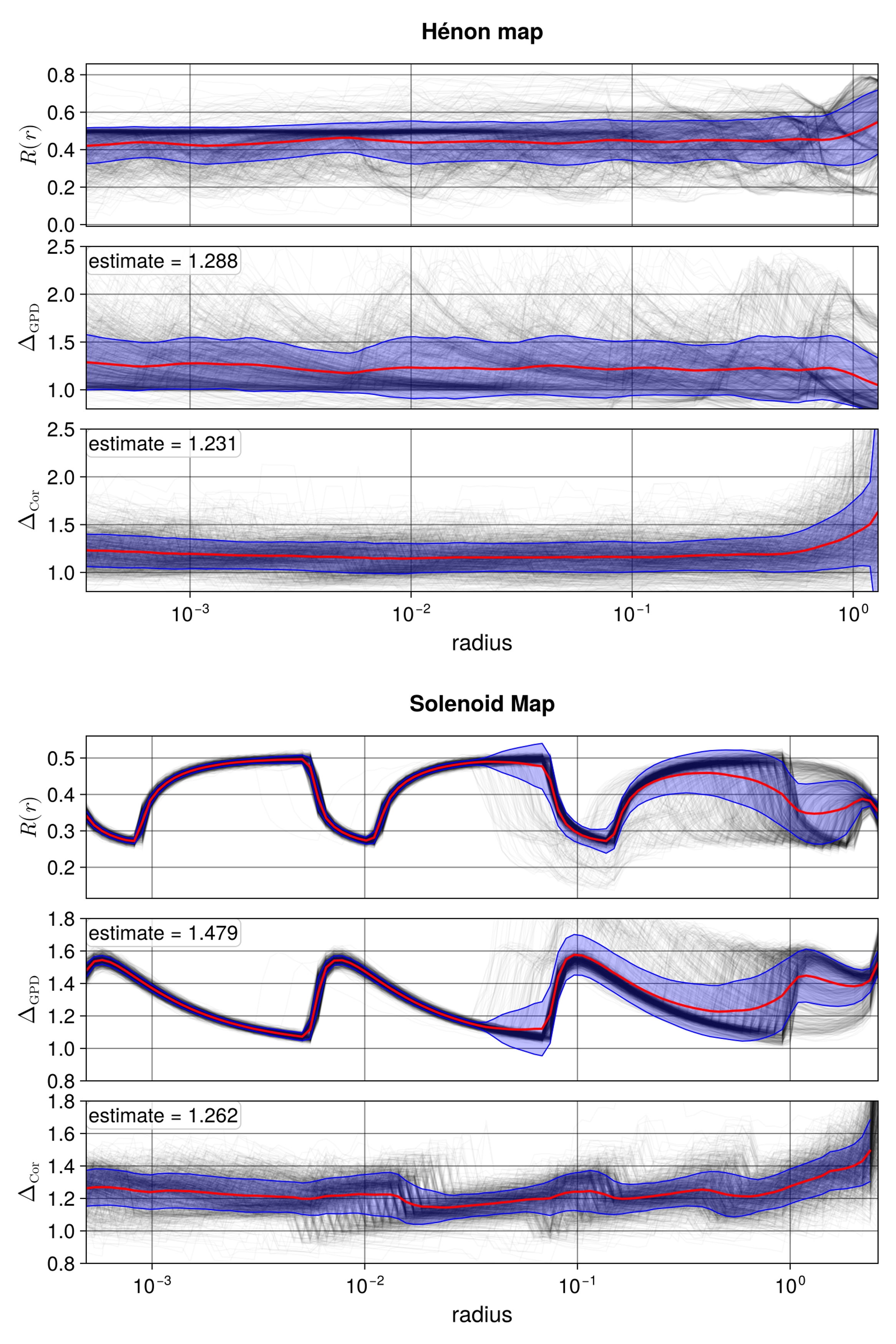}
    \caption{Same kind of figure as before, but displaying the quantities for 1000 different randomly selected points in black, their ensemble average in red and its standard deviation as a range in blue. The top panel shows the H\'{e}non map and the bottom panel the Solenoid map.}
    \label{fig:discrete_maps}
\end{figure*}

Another use of this algorithm is to compute the information dimension of the attractor \cite{lucarini2016extremes,datseris2023estimating} by averaging the local dimension estimates along a trajectory. This typically gives values around $1.35$ for the H\'{e}non attractor, which is only one decimal point away for the value $1.26\pm 0.02$ estimated by other authors \cite{grassberger1983generalized}. This is possibly a reason why the this system has been used to illustrate the EBD method in some texts \cite{lucarini2014towards,lucarini2016extremes}, despite the lack of regular variation. The top panel of Figure~\ref{fig:discrete_maps} shows the result of zooming and averaging over 1000 randomly selected reference points. The quotient of the measure of the two balls described by Eq.~\eqref{eq:regularly_varying_quantity} oscillates for any given point, at different values of the radius. There is a clear cluster around the value $2^{-1}$, and corresponding $d_\zeta=1$ meaning that the set of points that lie inside the ball looks like a line at that resolution. While the average value of the dimension (approximately $1.255$) keeps reasonably flat and very close to the value estimated by other sources \cite{grassberger1983generalized}, we have to remind the reader that this value however is not the output of applying the algorithm to the process given in Eq.~\eqref{eq:algorithm1process}, and actually used in applications, but the output of applying it to the process in Eq.~\eqref{eq:hittingprocess}, since we are not averaging over estimates done along an orbit, but rather averaging over the estimates performed on independently chosen points, each of them getting its own orbit, and drawing only one estimate from each orbit. This does indeed make a difference, since the output of the algorithm along one orbit gives $1.361$, while the average over independent points give $1.255$.

\subsection{Axiom A Systems: The Solenoid}

As a last example for this section we have the Solenoid map. This map is and example of a uniformly hyperbolic system, also called sometimes axiom A. These systems are usually presented as models for more complicated systems, and have been presented in the literature as a valid setting to apply this algorithm\cite{lucarini2016extremes,lucarini2014towards}. In \cite{pons2023statistical} it is stated that the exceedances of the process \eqref{eq:hittingprocess} around a point in the attractor of an Axiom A system converge to an exponential distribution with parameter equal to the inverse of the local dimension around that point. Here we produce a counterexample.

The solenoid map is defined as a transformation on a torus. Following the construction in \cite{falconer2007fractal}, let $V = \{(x,y): x^2+y^2 \leq 1\}$ be the unit disk, and $W$ the unit circle. Then we can parametrize a solid torus as 
\begin{align*}
    \mathbb{T}_S= \{(\varphi,v)\in W\times V: 0\leq\varphi<2\pi, |v|< 1\}
\end{align*} 
and the Solenoid map can be defined as 
\begin{align}\label{eq:Solenoid_rule}
    f_S(\varphi,v) = (2\varphi (\text{mod }2\pi), av + \overline{\varphi}/2)
\end{align}
where $\overline{\varphi}$ is the unit vector in $V$ at an angle $\varphi$ with the $x$ axis, and $a$ is a number $0<a<1/4$. It can be shown that this system has an attracting invariant set which can be defined as 
\begin{align}\label{eq:Solenoid_attractor}
    A_S = \bigcap_{i=0}^\infty f_S^i(\mathbb{T}_S)
\end{align}

and whose Hausdorff dimension is equal to 
\begin{align}\label{eq:Solenoid_dimension}
    \dim_H(A_S) = 1 - \frac{\log 2}{\log a},
\end{align}
see \cite{falconer2007fractal}. Note that for Axiom A systems the local dimension is constant $\mu$ almost everywhere and equal to the Hausdorff dimension, the Box-Counting dimension, information dimension, and most of the frequently used notions of dimension \cite{young1982dimension,barreira2011dimension,barreira2001hausdorff}.
An approximation of the measure of a ball according to the invariant measure of the solenoid $\mu_S$ can be derived by introducing some simplifications. The derivation is explained in detail in appendix \ref{ap:estimation_solenoid_measure}. Using this approximation, the measure of a ball of radius $r$ centred around a point $\varphi,v$, behaves approximately as a power law $\mu_S(B_r(\varphi,v))\approx r^{1-\log 2/\log a}$, but with a series of kinks spaced according to a geometrical series $r_n= h_a a^n$ where $h_a$ is a constant which may depend in the parameter $a$. This is consistent with the geometrical phenomenon described in figure~\ref{fig:HenonZooms}, since clusters of points in the Poincar\'{e} section $A_S\cap P_\varphi$ where $P_\varphi$ is a semi-plane at angle $\varphi$ from the centre, are separated by powers of $a$.


In the bottom panel of figure \ref{fig:discrete_maps} we the same plot as before with the H\'{e}non map but performed for the Solenoid map, with parameter $a=0.076$. The quotient of the measure of the balls becomes independent of the chosen reference point due to the uniformity of the attracting set $A_S$. The result is an oscillation which synchronizes for any reference point, and the same for the estimated local dimension. As a result, there is no convergence for the algorithm even on average. 
The approximation follows closely the numerical average of the quantity \eqref{eq:regularly_varying_quantity} computed over all the points, and allows us to extend it with much lower computational effort. This approximation can yield good values for very small radii, and as a verification, estimating the local dimension as the slope of the plot $\log \mu_S(B_r(\varphi,v))$ vs. $\log r$ over 16 orders of magnitude gives $1.2695\dots$, which is very close to the exact value $1-\log 2/\log 0.076 = 1.2689\dots$, which is an accuracy hard to reach when estimating dimensions of fractal sets.

\section{Continuous Systems}

\subsection{Extremal index estimates}

A continuous time notion of the extremal index exists \cite{fasen2009extremes}, however, previous work\cite{alberti2023dynamical,faranda2017dynamical,faranda2017dynamicalB,faranda2020diagnosing,hochman2019new,messori2017dynamical,rodrigues2018dynamical} that analyses data through this methodology compute the discrete time version of the extremal index under the assumption that the data is a discrete sampling of a continuous process. The papers use the index based estimator given by \cite{suveges2007likelihood}, which infers the value of $\theta$ through the indexes of the extremes as in a discrete time process $X_1,X_2,\dots$, and thus there is an ambiguity in the method used to discretize the continuous time process into a discrete time one. While this is unambiguous if one uses the return map to a Poincar\'{e} surface or a suspension flow, problems can arise if one uses a numerical scheme as temporal discretization. Intuitively, the extremal index describes the inverse average size of clusters, and the estimator infers this from the index length between exceedances. Clustered exceedances correspond to a lower value of $\theta$, meanwhile longer gaps between them will increase it.

The translation from continuous time to indexes of a discrete process depends heavily on the sampling frequency of the data and the length of the data series, and at the same time it has an interplay with the threshold choice. Denser sampling will make clusters bigger, but also make the gaps between clusters bigger in the index sense. Longer trajectories may decrease the sizes of the clusters, but they will introduce more clusters and more gaps of unknown length. Increasing the threshold will make the clusters smaller and the gaps bigger, therefore increasing the value of $\theta$.

\begin{figure*}
    \centering
    \includegraphics[width=\linewidth]{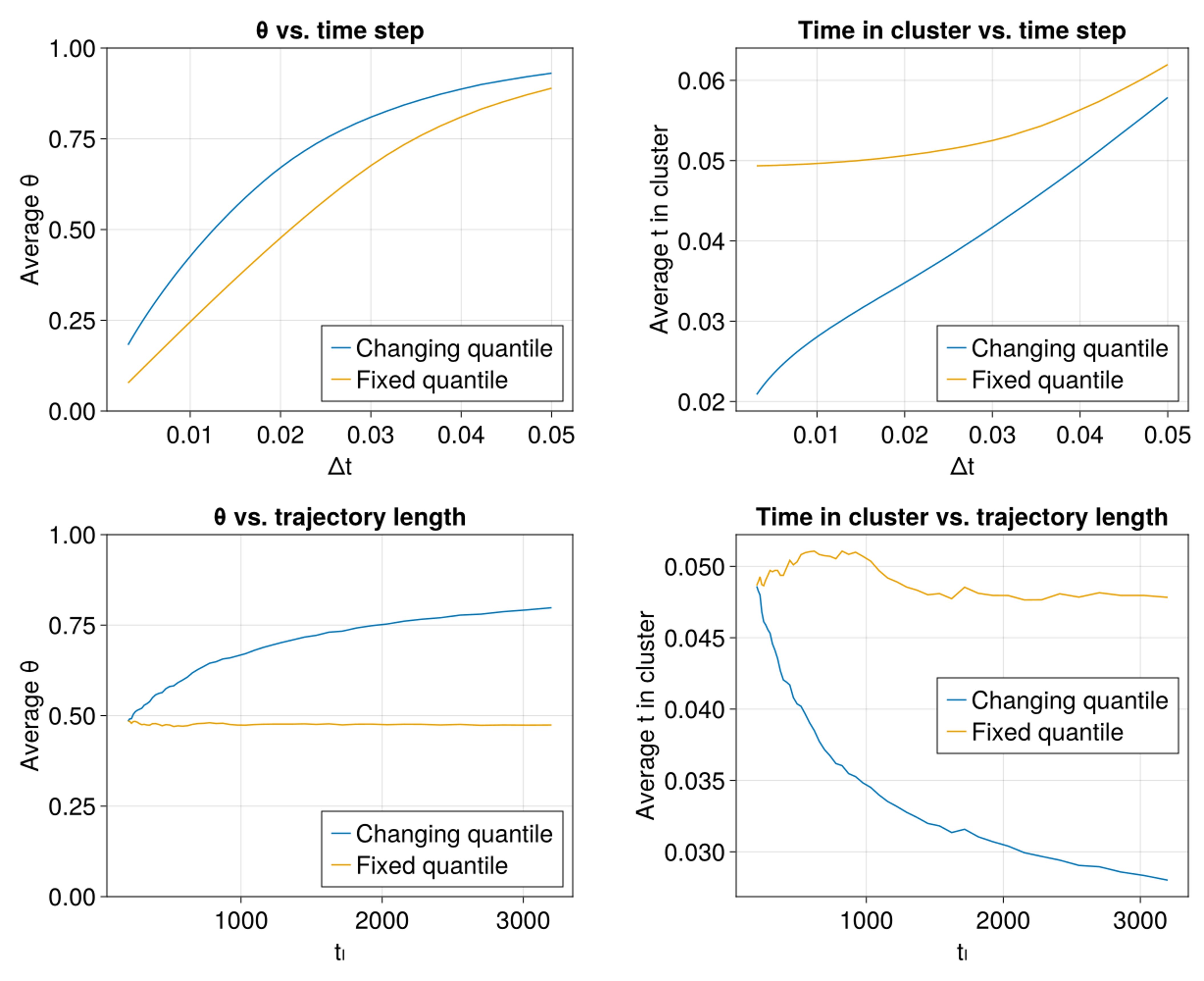}
    \caption{On the left, averaged extremal index computed over the points in a trajectory using the index estimator. On the top panel it is computed along trajectories of the same temporal length but with different sampling times, the bottom panel shows the effect of changing the trajectory temporal length with a fixed time step of the integration. On the right, averaged time spent in a cluster computed over the points in a trajectory using the extremal index. On the top panel it is computed along trajectories of the same temporal length but with different sampling times, the bottom panel shows the effect of changing the trajectory temporal length with a fixed time step of the integration.}
    \label{fig:extremal_index}
\end{figure*}

To study the effects of varying this parameters, we have made a simple experiment with the Lorenz 63 model, the result of which can be seen in figure ~\ref{fig:extremal_index}. On the top left panel the extremal index is computed and averaged over the process \eqref{eq:algorithm1process} using points on a trajectory for different values of the time step $\Delta t$ of the integration and a fixed temporal length of $t_l = 1000$ time units. Two curves have been computed, one with a fixed quantile $q=0.99$ and one with a varying quantile $q = 1-1/\sqrt{t_l/\Delta t}$, which is the quantile that corresponds at taking as extremes approximately the $\sqrt{N}$ biggest values, where $N\approx t_l/\Delta t$ is the total number of data in the process. For both curves the extremal index decreases strongly, given that $0\leq\theta\leq 1$, as the time step decreases. The bottom left panel shows a similar experiment but with a fixed time step $\Delta t\approx 0.0198$ and for different the temporal lengths of the trajectory. The quantiles are the same as in the top left panel. Now there is a difference in the behaviour of the two curves. While the extremal index estimation does not change with length when the quantile stays fixed, the extremal index grows with length when the quantile also increases. 

This problem has been acknowledged in most of the papers dealing with applications\cite{faranda2020diagnosing,hochman2019new,faranda2017dynamicalB,rodrigues2018dynamical}. A proposed solution is to renormalize by the time step to obtain a coherent temporal duration independent of the discretization, even if the size of the cluster varies in the index sense. Figure~\ref{fig:extremal_index} shows that it is not the case. By the interpretation that the extremal index is the inverse of the average size of clusters of exceedances we can "invert" the discretization back into continuous time. Calling $i_c$ the cluster length counted in indexes of the discrete time process, then the temporal length of clusters $t_c$ in the continuous time is approximately

\begin{align}
    \theta = \frac{1}{\langle i_c\rangle} \Rightarrow \langle t_c \rangle \approx \langle i_c\rangle \Delta t = \frac{\Delta t }{\theta}.
\end{align}
This quantity, averaged over the trajectory, is shown in the right panels of Figure~\ref{fig:extremal_index}. In the top right panel we see that the time duration of the clusters decreases with denser sampling, although this time the behaviour is more affected by the quantile choice. When keeping the quantile fixed the time spent in the cluster decreases, but seems to stabilize as $\Delta t\to 0$. When the quantile changes with the amount of data, the average time in clusters decreases almost linearly.
This shows that in practice we cannot compare results from index based estimators over time series sampled with different frequency. The bottom right panel shows how the average time in cluster changes with the length of the trajectory. If the quantile is fixed, the average time spent in clusters remains approximately constant, but if the quantile increases the average time spent in clusters decreases with the length of the trajectory. 

The choice of quantile has a significant impact in the results. Generally there is a trade off between the quantile and the number of data points, while higher quantiles deliver a better approximation to the asymptotic limit that is the objective of the estimation, they also result in many points being discarded, and thus smaller sample sizes from which to perform statistical inference.
One of the reasons why the changing quantile chosen leaves $\sqrt{N}$ points as extremes is that increasing $N$ simultaneously increases the quantile and the size of the sample, and thus should give a better approximation in both fronts, but this can only be done with synthetic data. In most applications\cite{faranda2020diagnosing,faranda2017dynamicalB,faranda2017dynamical,hochman2019new,messori2017dynamical} the quantile is fixed and equal to $q=0.98$ although others use other fixed value \cite{alberti2023dynamical,rodrigues2018dynamical}.

An unambiguous discretization method such as Poincar\'{e} return map assigns an extremal index of $\theta = 1$ to a.e. point for the Lorenz 63 system attractor\cite{zhang2015borel}.

\subsection{Lorenz 63}

\begin{figure}
    \centering
    \includegraphics[width=\linewidth]{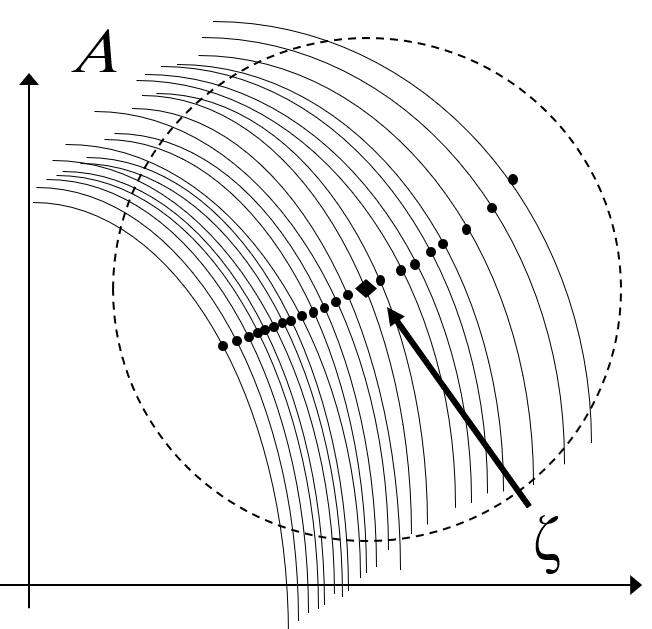}
    \caption{Schematic representation of computational method used to study continuous time dynamical systems. The reference point $\zeta$ is represented by a diamond shaped point and circular points are the closest recurrences to it of each branch of the attractor $A$ that enters the ball represented with a dashed line.}
    \label{fig:SchemaContinuousSystems}
\end{figure}

So far, we are iterating the dynamical system until we return to a small neighbourhood close to a reference point. In continuous time systems, this comes with a computational difficulty. Numerical schemes introduce a source of error by producing points in the trajectory which may be too far from the reference point even if the continuous time trajectory goes close enough to the reference point. As the radius of the ball decreases, the numerical approximation may even miss to hit the ball and omit a true recurrence. These problems can be avoided by performing an interpolation and solving a minimization problem that produces the only the closest recurrence of the trajectory within the ball. This means that for every trajectory that enters the ball, only one point is kept. Thus the set of points that remain in the ball has a dimension which is the dimension of the intersection between the attractor and ball minus 1. Figure \ref{fig:SchemaContinuousSystems} shows a schema of an attractor intersecting a ball centred around a reference point and the points that remain inside it after the optimization procedure takes place. We use this method to analyse the regularly varying properties of the continuous time systems that appear in this section from now on.

\begin{figure*}
    \centering
    \includegraphics[width=0.8\linewidth]{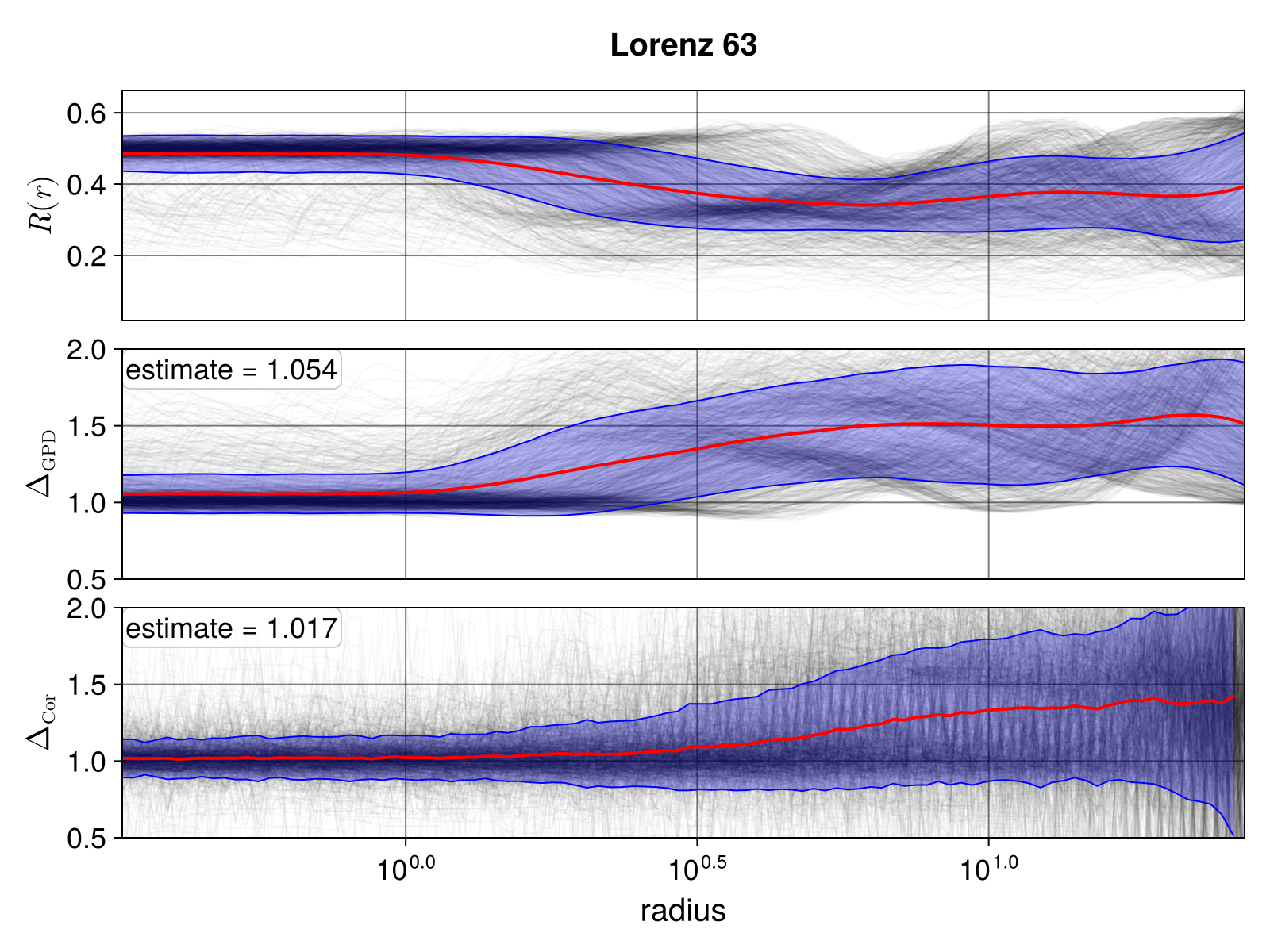}
    \caption{This figure shows the quantity $R(r)$ in the top panel, the corresponding GPD estimation of the local dimension in the middle panel and the correlation estimation of the local dimension in the bottom panel for a 1000 different points randomly selected along the attractor.}
    \label{fig:Lorenz63ManyPoints}
\end{figure*}

The Lorenz 63 system is described by the equations
\begin{align}
    \begin{split}
    &\Dot{x_1} = \sigma(x_2-x_1) \\
    &\Dot{x_2} = x_1(\rho-x_3)-x_2\\
    &\Dot{x_3} = x_1x_2-\beta x_3
    \end{split}
\end{align}
with the parameter values $\sigma =10$, $\rho = 28$ and $\beta = 8/3$, and it has also been used to illustrate the EBD method\cite{faranda2017dynamical,rodrigues2018dynamical}. Figure~\ref{fig:Lorenz63ManyPoints} shows the analysis as performed with the discrete systems but with the interpolation scheme. As the radius decreases the estimates of the measure ratio cluster around $0.5$ with some downward excursions. Similarly, the estimates of the dimension cluster around 1, but have oscillations that push them upwards for some range of scales. This phenomenon is similar to what was observed for the H\'{e}non map, indicating that the attractor is locally very flat and behaves like a surface for most points at any given sufficiently small scale, but not for any point through all the scales. The average of the local dimensions stabilizes around the value $2.054$, very close to the value $2.06\pm 0.01$ reported in \cite{grassberger1983characterization}.

Again, this means that the invariant measure of the system is not regularly varying and thus does not fulfil the conditions needed for the EBD algorithm to work, even if the average result is very close. The individual estimations of the local dimension around a point are mostly very close to 2, the value corresponding to a surface with the rest of the points incrementing slightly the value. It is remarkable how stable this configuration is across scales once a sufficiently small scale is reached around $r\simeq 1$. The results are similar for the correlation estimate, which is unstable, and also gives a value close to 2. Similar results are obtained for other very flat attractors such as the one arising in Rössler system.

\subsection{Lorenz 96}

\begin{figure*}
    \centering
    \includegraphics[height=0.9\textheight]{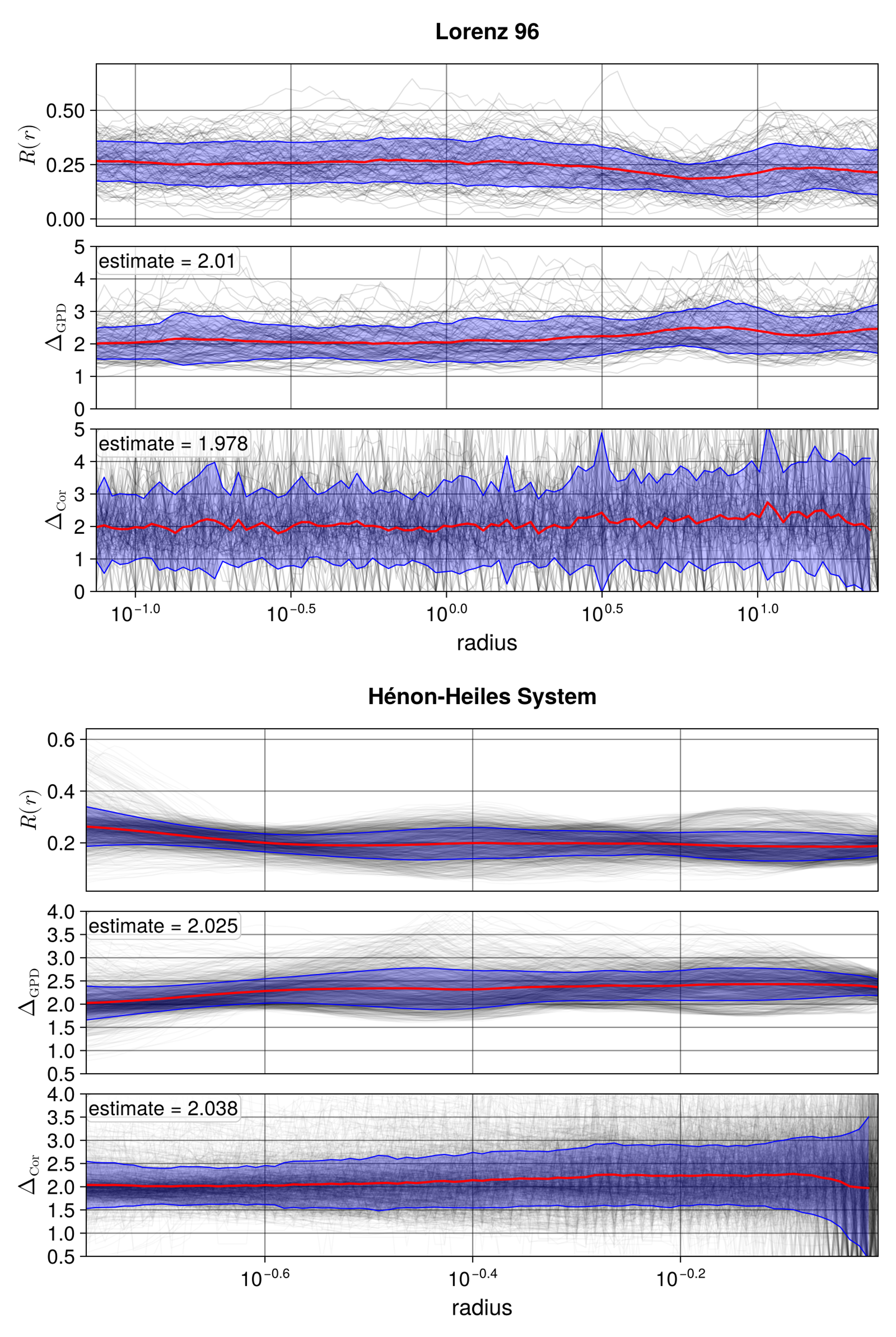}
    \caption{Same as the rest but with the 4D Lorenz 96 model on the top panel and the H\'{e}non-Heiles system in the bottom panel.}
    \label{fig:continuous_time_regularly}
\end{figure*}

An example of a less flat attractor is given by the Lorenz 96 model. This is given by the equations 
\begin{align}
    \Dot{x_j} = x_{j-1}(x_{j+1}-x_{j-2}) - x_j + F
\end{align}
with $1\leq j \leq n$, and assuming that $x_{n+1} = x_1$, $x_0=x_n$ and $x_{-1}=x_{n-1}$. The top panel of Figure~\ref{fig:continuous_time_regularly} shows the results for this system with the parameters $n=4$ and $F=32$, which results in chaotic time evolution. The quantities plotted do not seem so converge or cluster in any particular way. A sufficiently small scale seems to be reached around $r=10^{0.5}$, after that the distribution of estimates remains stable. The range of the estimations remains very broad and with a big dispersion. Our estimation of the information dimension of the attractor over a trajectory with $10^7$ points is $2.834$, which differs significantly from the average answer of $3.010$ of the run, and from the actual result of applying the EBD algorithm taking the reference points in a trajectory, which gives $3.023$. The results are similar for the correlation estimate, which is very unstable for this system. 

\subsection{H\'{e}non-Heiles system}

The H\'{e}non-Heiles system is a Hamiltonian system that was introduced to model the orbit of stars around galactic centres, and its orbits also display a fat fractal behaviour \cite{barrio2008fractal}. The equations of motion of the the H\'{e}non-Heiles system are given by 

\begin{align}
    \begin{split}
    &
    \Dot{x} = p_x \\
    &\Dot{p_x} = -x-2xy\\
    &\Dot{y} = p_y\\
    &\Dot{p_y} = -y - (x^2-y^2)
    \end{split}.
\end{align}
The initial condition provided to this system determines the total energy available to the orbits, for our simulations we have used far apart initial conditions belonging to the trajectory with initial condition $(0.0, -0.25, 0.42, 0.0)$.


The results for this system are not very conclusive. A problem that arises is that the numerical integration produces round up errors that do not conserve the energy of the system on the long run. This means in our case that the reference point is in a manifold from which the dynamics drift away. We would expect to see convergence in this case, and indeed the top graph and the middle graph of the bottom panel of figure~\ref{fig:continuous_time_regularly} could be converging around $r\approx10^{-0.6}$, but start diverging significant for smaller radii. ottom graph gives the correlation estimation of the local dimension, which seems to be quite unstable, but averages close to 3.


\section{Conclusions}

The main question that this paper aims to answer deals with the abundance of systems for which the GPD estimators of the local dimension are adequate. In the theoretical papers that developed the EBD method an abstract proof is provided which references the mathematical properties that need to hold in order for the EBD algorithm to work. However, the properties are just assumed to hold for large classes of dynamical systems, and are not checked either in the examples provided or in the data to which it is applied. Moreover, in the papers that deal with applications to data there are several steps requiring further mathematical justification. This corresponds to whether an  extremal index should be included in the GPD distribution function, on the performance of dimension estimation schemes when changing the process given in Eq.~\eqref{eq:hittingprocess} by the process given in Eq.~\eqref{eq:algorithm1process}, or setting the tail parameter $\xi\neq 0$ when the system is not uniformly hyperbolic \cite{pons2023statistical}. 

One of the focuses of this paper is on the regular variation property because it is probably the most uncommon one of the required properties listed above. It is also the least studied theoretically, when compared to the prolific mathematical theory linking local dimension, ergodicity, return times and other characteristics of chaotic systems. Also the focus is in the behaviour of the estimation for almost every reference point, since since fixed points and periodic orbits have an important effect in the method, but are usually a measure zero set which is unlikely to be observed in applications.

Similarly, problems arise when considering other dynamical indicator of persistence, namely the extremal index. For discrete systems, the extremal index is 1 $\mu$-a.e, except in a zero measure set of special points. For continuous systems sampled at a fixed time intervals, the extremal index estimation depends heavily on parameters such as the sampling frequency and the threshold choice, which one may not even have control over when working with data.

We have shown here some experiments that we believe are particularly meaningful and help to the understanding of the relation between the geometry of the attractor and the regular variation property. Among the systems tested but not shown are chaotic discrete maps, a quasiperiodic discrete map with a strange attractor \cite{grebogi1984strange}, intermittent discrete maps, continuous chaotic systems and a non-autonomous quasiperiodic continuous system. From all this experiments transpires that the geometry of the attractor is the determining factor that leads to the regular variation property holding or not, and thus the EBD algorithm working. The algorithm showed good performance for the quasiperiodic examples, except for the case that had a strange attractor. Other differences such as chaotic/not chaotic, discrete/continuous, autonomous/non-autonomous did not made a difference. Particularly, we could not find any regularly varying measure in systems such that the asymptotic invariant measure of an orbit was singular with respect to Lebesgue, thus in any system having a dimension different from an integer.

The implications of our work for the application of the EBD algorithm is in the eye of the beholder.
Primarily, our goal is to highlight that further checks are needed to obtain confidence in applying the EBD algorithm to data. The method appears in many cases to output reasonable, or expected, answers, but such answers may be just an artifact coming from free choices such as parameters (e.g. the threshold $q$) or arbitrary quantities (e.g. the length of the data series or the sampling frequency).
Even in cases considered in this paper where the criteria for the EBD method are not satisfied, the estimates that are shown in this paper are quite close to the theoretical or expected results. One possible reason for this is the low dimensionality of the systems used. It is computationally expensive to check higher dimensional systems, where the differences between the estimation and the true result could be greater, since the number of iterates needed to decrease the radius becomes large.


\appendix 
\section{Correlation-based estimation for the local dimension} \label{ap:correlation_dimension_explanation}

The correlation method to estimate the local dimension is based on the decay of the a quantity called the correlation sum, given by
\begin{align}\label{eq:correlation_sum_order_2}
    S_C(\zeta,\epsilon) = \frac{1}{N}\sum_{i=1}^N \mathbbm{1}_{\{||\boldsymbol{x}_i-\zeta||<\epsilon\}}(\boldsymbol{x}_i)
\end{align}
where $\mathbbm{1}_A(\boldsymbol{x})$ is the indicator function of the set $A$, and thus its value is $1$ if $\boldsymbol{x}\in A$ and $0$ otherwise. The scaling of this quantity with its argument $\epsilon$ for small $\epsilon>0$ is equivalent to the scaling of the measure around the point $\zeta$, and thus one has
\begin{align}
    \lim_{\epsilon\to 0}\frac{\log S_C(\zeta,\epsilon)}{\log \epsilon}\to \Delta^l_\zeta.
\end{align}

To estimate the local dimension, we compute the correlation sum around the reference point $\zeta$ and make a log-log plot of the correlation sum against $\epsilon$. Then compute the largest linear region in the plot and fit a linear regression to it. The slope of the linear regression is the estimate of the local dimension around $\zeta$. The correlation method can be presented with much more generality by introducing an order $s\in\mathbb{N}$, but it is sufficient for our purposes to consider only one order $s=2$, for which the correlation sum takes the simple form shown in \eqref{eq:correlation_sum_order_2}. For a more in depth explanation of the correlation method and its comparison with other estimation procedures for fractal dimensions see\cite{datseris2023estimating}.

\section{Estimating the invariant measure of the Solenoid Map}\label{ap:estimation_solenoid_measure}

The proof of formula \eqref{eq:Solenoid_dimension} in \cite{falconer2007fractal} is made by finding a suitable covering by balls of an iterate $f_S^k(\mathbb{T}_S)$. This idea can be modified to estimate the volume of a small ball intersecting the Solenoid attractor. The solenoid attractor embedded in $\mathbb{R}^3$ looks like a line coiled infinitely many times within a torus, hence the name, and the intersection of the solenoid attractor with a small ball centered in a point in the attractor has the structure of a series of branches of variable length, slightly curved, one of which goes through the center. We do the following simplifications:
\begin{itemize}
    \item The measure is uniform along the 
    length of the branches.
    \item Since we are interested in the measure of balls of a very small radius, we assume that the branches within the ball are approximately straight.
    \item All branches are of total length $2\pi$ when they go around the torus.
\end{itemize}

This simplifications imply that we approximate the invariant measure of the Solenoid map with a measure on a Cantor comb, formed by the product $[0,2\pi]\times P_{\varphi}$ where $P_\varphi$ is the Poincar\'{e} section of the Solenoid map of a semi-plane at angle $\varphi$, effectively, we are computing how much of the length of series of circles disposed to form a Poincar\'{e} section in one semi plane is covered by ball of small radius.  To obtain an approximation of the Poincar\'{e} section fix a $k\in \mathbb{N}$, and iterate forward the points that will land in $P_\varphi$ after $k$ iterations. Since the Solenoid map acts as the doubling map in the angle, the image of the torus intersects 2 times with $P_\varphi$, and the $k$th image of the torus will intersect $2^k$ times with it. Writing as a sequence $(\varphi_n,v_n)=f_S(\varphi_{n-1},v_{n-1})$ the angle $\varphi_k$ iterated backwards gives
\begin{align}\label{eq:angle_values}
    \begin{aligned}
    &\varphi_{k-1}^{\Gamma_{k}} = \varphi_k/2 + a_k \pi\\
    &\varphi_{k-2}^{\Gamma_{k-1,k}} = (\varphi_k/2 + a_k\pi)/2 + a_{k-1}\pi\\
    &\cvdots\\
    & \varphi_0^{\Gamma_{1,k}} = \varphi_k/2^k + \sum_{i = 1}^k \frac{a_i\pi}{2^{i-1}}
\end{aligned}
\end{align}
where $\Gamma_{i,j}=\{a_i,\dots,a_j\}\in \{0,1\}^{j-i}$ is a sequence and the values $a_i$ are either $0$ or $1$, and $\Gamma_{i}=\Gamma_{i,i}=\{a_i\}$. Note that there are $2^{j-i}$ possible sets $\Gamma_{i,j}$ and thus also $2^{j-i}$ possible values for the angle $\varphi_{i-1}^{\Gamma_{i,j}}$, corresponding to a particular realisation of the sequence $\Gamma_{i,j}$. In this setting,
we are interested in looking at the Poincar\'{e} section $P_{\phi_k}$, and the angles $\varphi_0^{\Gamma_{1,k}}$ are exactly the angles that will land in $P_{\phi_k}$ after $k$ iterations of the map. The contribution of the angle $\varphi_k$ is exponentially small, and thus for high $k$ the angles $\varphi_0^{\Gamma_{1,k}}$ are very weakly influenced by the Poincar\'{e} section that we are looking at.

To find out the geometrical structure of the points in $P_{\varphi_k}$ now we compute the dynamic law but forward this time in the $v$ coordinate. Given an initial condition $v_0$
\begin{align}\label{eq:circle_point_values}
    \begin{aligned}
    &v_1^{\Gamma_{1,k}} = av_0 + \overline{\varphi}^{\Gamma_{1,k}}_0/2\\
    &v_2^{\Gamma_{1,k}} = a(av_0+ \overline{\varphi}^{\Gamma_{1,k}}_0/2) + \overline{\varphi}_{1}^{\Gamma_{2,k}}/2\\
    &\cvdots\\
    & v_k^{\Gamma_{1,k}} = a^kv_0 + \sum_{i = 1}^{k} \frac{\overline{\varphi}^{\Gamma_{i,k}}_{i-1} a^{k-i}}{2}
\end{aligned}
\end{align}
where the unitary vectors $\overline{\varphi}_{j-1}^{\Gamma_{j,k}}$ are defined as 
\begin{align*}
    \overline{\varphi}_{j-1}^{\Gamma_{j,k}}= \left( \cos  \left(\varphi_k/2^{k-j} + \sum_{i = j}^{k-1} \frac{a_i\pi}{2^{i-1}}\right) ,
     \sin\left( \varphi_k/2^{k-j} + \sum_{i = j}^{k-1} \frac{a_i\pi}{2^{i-1}} \right) \right).
\end{align*}

Equation \eqref{eq:circle_point_values} gives an approximation of where the set $f_S^k(\mathbb{T}_S)$ intersects $P_{\varphi_k}$, since $v_0\in V$, we have $v_k^{\Gamma_{1,k}}\in f_S^k(\mathbb{T}_S)$, and thus its distance to the nearest point of $A_S\cap P_\varphi$ is bounded by $2a^k\leq 2/4^k$, corresponding to the maximal diameter of a connected component of $f_S^k(\mathbb{T}_S)\cap P_{\varphi_k}$. 

Now using the simplifications stated earlier the measure of a ball centered on a point $(\varphi_k,v_k^{\Gamma^*})$ for some $\Gamma^*\in \{0,1\}^k$ and $r>0$ can be estimated as 
\begin{align}\label{eq:Solenoid_measure}
    \mu_S(B_r(\varphi_k,v_k^{\Gamma^*}))= \frac{1}{2^k\pi}\sum_{\Gamma_{1,k}}\sqrt{r^2-||v_k^{\Gamma^*}-v_k^{\Gamma_{1,k}}||}
\end{align}
with the convention that if $r^2<||v_k^{\Gamma^*}-v_k^{\Gamma_{1,k}}||$ for a sequence $\Gamma_{1,k}$ that term does not contribute to the sum. The approximation given by formula \eqref{eq:Solenoid_measure} is not very useful analytically since it is cumbersome to work with (having an exponential amount of terms, having to compute the distances between all points to know if they affect the sum or not, not knowing the speed of convergence in $k$ or proving that it does converge...) and it is based in simplifications which do not allow to quantify the error easily. The good news is that numerically it seems to work very accurately and this formula also sheds some light in the behaviour of the actual measure.

One thing that becomes clear on the derivation of the approximation of the measure is that the point chosen to do the approximation does not matter much. This is, in formula \eqref{eq:angle_values} we see that the initial angle at which we take the Poincar\'{e} section is exponentially dampened, so the angles will be almost the same for any Poincar\'{e} section taken. The initial condition within the disk $v_0\in V$ is also exponentially dampened because $a<1/4$, so again it does not contribute much for high $k$. This is consistent with the uniform hyperbolicity of the system, which implies that the invariant set should have locally the same structure everywhere. 

\section{Data availability statement}

The code and the data used to generate the plots in this paper is fully available on Github\cite{CodeBase}.

\section*{Acknowledgements}
This research is supported by the
European Union’s Horizon 2020 research and innovation programme under
the Marie Sklodowska-Curie Actions grant No.956170 CriticalEarth.

\bibliography{REFERENCES}

\end{document}